# EXACT SOLUTIONS TO THE NAVIER-STOKES EQUATION FOR AN INCOMPRESSIBLE FLOW FROM THE INTERPRETATION OF THE SCHRÖDINGER WAVE FUNCTION

Vladimir V. KULISH* & José L. LAGE**

*School of Mechanical & Aerospace Engineering, Nanyang Technological University, Singapore 639798
** Lyle School of Engineering, Department of Mechanical Engineering, Southern Methodist University, Dallas TX 75275-0337

January 17, 2013

**Abstract**

The existence of the velocity potential is a direct consequence from the derivation of the continuity equation from the Schrödinger equation. This implies that the Cole-Hopf transformation is applicable to the Navier-Stokes equation for an incompressible flow and allows reducing the Navier-Stokes equation to the Einstein-Kolmogorov equation, in which the reaction term depends on the pressure. The solution to the resulting equation, and to the Navier-Stokes equation as well, can then be written in terms of the Green's function of the heat equation and is given in the form of an integral mapping. Such a form of the solution makes bifurcation period doubling possible, i.e. solutions to transition and turbulent flow regimes in spite of the existence of the velocity potential.

## 1. Introduction

The Navier-Stokes equation, together with the continuity equation, is the governing equation of fluid motion valid as long as the continuum hypothesis holds. For a wide variety of flow phenomena encountered in science and engineering, the assumption of constant properties is adequate. However, even when simplified by this assumption the Navier-Stokes equation still presents a challenge, namely the nonlinear convective term, which makes the solving procedure very cumbersome, if possible at all.

It is commonly believed that the Navier-Stokes equation indeed governs *any* fluid motion, including turbulence [F]. If it is so, no additional hypothesis should be made in order to describe turbulence by means of the Navier-Stokes equation. However, since no general method for solving the Navier-Stokes equation exists nowadays, it has become usual to resort to various models that often are inconsistent with the phenomenon in question. Thus, for example, the turbulence concept of eddy viscosity emerges independent of any physical properties of the fluid. The advantage of splitting the velocity components into time averaged and fluctuating parts is doubtful, for it leads to the closure problem of the Navier-Stokes equation. Besides, it is not completely clear how the averaging should take place in reality, for depending on the scale of the time step used for averaging the result can differ significantly.

In our previous work, we demonstrated that the Navier-Stokes equation reduces to the Einstein-Kolmogorov (reaction-diffusion) equation under the assumption of the existence of narrow probability wave packets [K-L].

In this work we demonstrate that, without any further assumptions, if the Schrödinger equation is valid to describe a physical state of a fluid particle, the continuity equation implies the existence of the velocity potential. This potential is nothing else but the phase of the corresponding wave function. Hence, the Cole-Hopf transformation [C], [H], in which the velocity function is expressed through the corresponding Schrödinger wave function, reduces the Navier-Stokes equation to a reaction-diffusion equation in the form of the Einstein-Kolmogorov equation. The solution to the resulting equation, and, by consequence, to the Navier-Stokes equation, can then be written in terms of the Green's function of the heat equation and is given in the form of an integral mapping. Such a form of the solution makes bifurcation period doubling possible, i.e. solutions for transition and turbulent flow regimes in spite of the existence of the velocity potential.

## 2. Existence of the velocity potential

The purpose of this section is to demonstrate that, as long as the continuity equation holds, the instantaneous velocity of the flow can *always* be represented as $\mathbf{v} = \nabla \widetilde{\theta}$, where $\widetilde{\theta}$ – the velocity potential – is nothing other than a quantity proportional to the phase of the wave function associated with the moving particle. The wave function is a complex function that satisfies the Schrödinger equation and is written as $\psi = Ae^{i\theta}$, where $A$ is the observable amplitude and $\theta$ is the unobservable phase.

The fact that the instantaneous velocity is proportional to the gradient of the phase of the corresponding wave function is well-known among the quantum mechanical community, [Hbrg], [Lgt]. In fact, this relationship can be established starting directly from the conservation equation, namely:

$$\frac{\partial w}{\partial t} + \operatorname{div} \mathbf{j} = 0$$

$$(2.1)$$

where $w$ denotes the energy density and $\mathbf{j}$ stands for the energy flux vector.

Equation (2.1) is the most general form of the conservation equation. Indeed, if one notes that $w$ can be interpreted as the density of the particles, then (2.1) is the equation describing conservation of the amount of particles; hence, $\mathbf{j}$ should be interpreted simply as the flux of particles. Furthermore, multiplying eqn. (2.1) by the mass of a particle $m$, one obtains

$$\frac{\partial \rho_m}{\partial t} + \text{div } \mathbf{j}_m = 0$$

(2.2)

where $\rho_m$ is the particle density and $\mathbf{j}_m$ is the density flux. Therefore, eqn. (2.2) is nothing other than the *continuity equation*. Even further, multiplying eqn. (2.1) by the electric charge of the particle, one obtains the equation describing conservation of the electric charge, i.e.,

$$\frac{\partial \rho_e}{\partial t} + \text{div } \mathbf{j}_e = 0$$

(2.3)

where $\rho_e$ is the density of the charge and $\mathbf{j}_e$ is the density of the electric current.

Recall now that $w = \psi\psi*$ by definition. Indeed, $\psi = Ae^{i\theta}$ and its conjugate $\psi* = Ae^{-i\theta}$ lead to $\psi\psi* = A^2 = w$, i.e., the probability density to find the particle in the state described by $\psi$. This is stated, for instance, in [B]. Since mass itself is a form of energy and the wave function $\psi$ represents the wave field of a particle of mass $m$, then $\psi\psi*$ may be thought of as an energy density associated with that mass. Accordingly, the product $\psi\psi* = |\psi|^2 = w$ is to be interpreted as a *probability density*.

Consequently, $\text{div } \mathbf{j}$ in (2.1) can be written as $\text{div } \mathbf{j} = iD\left(\psi^*\nabla^2\psi - \psi\nabla^2\psi^*\right)$ and $\frac{\partial w}{\partial t} = \left(\psi^* \frac{\partial \psi}{\partial t} + \psi \frac{\partial \psi^*}{\partial t}\right)$, where $D$ is a proportionality constant.

Hence, written in terms of the wave function, the conservation equation (2.1) becomes

$$\left(\psi^* \frac{\partial \psi}{\partial t} + \psi \frac{\partial \psi^*}{\partial t}\right) = iD\left(\psi^*\nabla^2\psi - \psi\nabla^2\psi^*\right)$$

(2.4)

Equation (2.4) can be split into two equations – one for the wave function and the other for its complex conjugate – namely:

$$\frac{\partial \psi}{\partial t} = iD\nabla^2\psi + \Phi(\mathbf{x},t)\psi$$

(2.5a)

and

$$-\frac{\partial \psi^*}{\partial t} = iD\nabla^2\psi^* + \Phi(\mathbf{x},t)\psi^*$$

(2.5b)

where $\Phi(\mathbf{x},t)$ is an arbitrary function.

Identifying $D$ as $\hbar/(2m)$, where $\hbar$ is the reduced Planck's constant, equations (2.5) become

$$i\hbar\frac{\partial\psi}{\partial t} = -\frac{\hbar^2}{2m}\nabla^2\psi + \Phi(\mathbf{x},t)\psi$$

(2.6a)

and

$$-i\hbar\frac{\partial\psi^*}{\partial t} = -\frac{\hbar^2}{2m}\nabla^2\psi^* + \Phi(\mathbf{x},t)\psi^*$$

(2.6b)

which are the Schrödinger equations for the wave function and its complex conjugate, respectively. In equations (2.6), therefore, $\Phi(\mathbf{x},t)$ should be understood as the field potential acting on the quantum system.

The flux in the conservation equation can now be written as

$$\mathbf{j} = -\frac{i\hbar}{2m}\left(\psi\nabla\psi^* - \psi^*\nabla\psi\right)$$

(2.7)

which, upon using $\psi = Ae^{i\theta}$ and $\psi^* = Ae^{-i\theta}$ can be written as

$$\mathbf{j} = A^2\frac{\hbar}{m}\nabla\theta$$

(2.8)

Now, since $A^2 = w$, then the quantity $(\hbar/m)\nabla\theta$ must be the velocity at the point $\mathbf{x}$.

From the analysis presented in this section, two important conclusions follow: firstly, provided the continuity equation holds, ***there always exists a function*** $\theta$ – the phase of the wave function – ***such that the instantaneous velocity*** $\mathbf{v}$ ***is directly proportional to the gradient of this function***, i.e., $\mathbf{v} = (\hbar/m)\nabla\theta$; and, secondly, ***the function used to perform the Cole-Hopf transformation***, $\mathbf{v} = -\dfrac{i\hbar}{m}\dfrac{\nabla\psi}{\psi}$, ***is the wave function in the Schrödinger equation that describes the quantum state of the velocity field and whose phase defines the velocity potential***.

Thus, velocity at a certain location $\mathbf{x}$ is proportional to the gradient of the phase of the corresponding wave function in the same location, i.e.,

$$\mathbf{v}(\mathbf{x}) = \nabla \widetilde{\theta} = \frac{\hbar}{m} \nabla \theta = -\frac{i\hbar}{m} \frac{\nabla \psi}{\psi}$$

(2.9)

where the value $(\hbar / m)\theta = \widetilde{\theta}$ is the *velocity potential*.

See [K-L] for a more elaborate procedure of obtaining the same result.

Hence, we demonstrated that *if the Schrödinger equation is valid to describe a physical state of a fluid particle, the continuity equation implies the existence of the velocity potential.* For a non-relativistic fluid particle (constant mass), this velocity potential is proportional to the phase of the corresponding probability wave function.

## 3. Applicability of the Cole-Hopf transformation to the Navier-Stokes equation for an incompressible fluid

It follows from the result, obtained in the previous section, that the velocity of the fluid particle is given in the form

$$\mathbf{v} = (\hbar / m)\nabla \theta$$

(3.1)

as long as the Schrödinger and continuity equations hold.

Now, because $(\hbar / m) = 2\nu$, where $\nu$ is the kinematic viscosity of the fluid in question, represents the same quantum state of the fluid particle, it follows that

$$\mathbf{v}(\mathbf{x}) = -2\nu \frac{\nabla \psi}{\psi}$$

(3.2)

The latter expression is known as the Cole-Hopf transformation [C], [H]. It was discovered by Hopf and Cole independently and is used to linearize the Burgers equation [L].

Upon substituting eqn. (3.2) into the Navier-Stokes equation written for an incompressible flow

$$\frac{\partial \mathbf{v}}{\partial t} + (\mathbf{v} \cdot \nabla)\mathbf{v} = -\frac{1}{\rho} \nabla p + \nu \nabla^2 \mathbf{v}$$

(3.3)

we obtain a reaction-diffusion equation, known as the Einstein-Kolmogorov equation [T-S]

$$\frac{\partial \psi}{\partial t} - \nu \nabla^2 \psi = \frac{\Delta p}{2\mu} \psi$$

(3.4)

where $\mu$ is the dynamic viscosity of the fluid in question, $p$ stands for pressure and $\Delta p$ denotes the pressure difference between the local pressure and a certain reference pressure, e.g. the ambient pressure. A detailed derivation of eqn. (3.4) from (3.3), by means of (3.2), is shown in [K-L].

Note that the corresponding transformations for the initial and boundary conditions are as follows

$$\psi(\mathbf{x},0) = C \exp\left( -i\frac{m}{\hbar} \int \mathbf{v}(\mathbf{x},0) d\mathbf{x} \right)$$

(3.5)

and

$$\psi = \exp\left( -\frac{1}{2\nu} \int \mathbf{v} \cdot d\mathbf{x} \right)$$

(3.6)

respectively.

## 4. Solution to the Navier-Stokes equation in the form of the Green's function of the corresponding heat equation

The global solution of eqn. (3.4) can now be written in terms of the Green's function of the non-homogeneous diffusion equation [Hman, p. 409]:

$$\psi(\mathbf{x},t) = \int_0^t \iiint \wp(\mathbf{x},t;\xi,\tau) \frac{\Delta p}{2\mu} \psi(\mathbf{x},t) d^3\xi d\tau + \iiint \wp(\mathbf{x},t;\xi,0) \psi(\xi,0) d^3\xi$$
$$+ \nu \int_0^t \oiint \left[ \wp(\mathbf{x},t;\xi,\tau) \nabla_\xi \psi - \psi(\xi,\tau) \nabla_\xi \wp(\mathbf{x},t;\xi,\tau) \right] \cdot \mathbf{n} dS d\tau$$

(4.1)

where $\wp(\mathbf{x}, t; \xi, \tau)$ is the Green's function of the diffusion equation for the domain of interest.

Observe that eqn. (4.1) has the form $\psi = Y(\psi)$, where Y is the mapping defined on the set of bounded continuous functions. Thus one can define a fixed-point, iterative process by $\psi_{m+1} = Y(\psi_m)$ with $\psi_0(\mathbf{x},t) = \iiint \wp(\mathbf{x},t;\xi,0)\psi(\xi,0)d^3\xi$ .

The initial approximation $\psi_0(x,t)$ is just the second term on the right side of eqn. (4.1), and it is the solution to the linear, homogeneous diffusion equation with the initial condition $\psi(x, 0) = \psi_0(x)$. Therefore, the solution of eqn. (3.3) can be found from (3.2), where $\psi(x, t)$ is determined by eqn. (4.1).

We note that the mapping form of the solution (4.1) makes bifurcation period doubling possible, i.e. solutions for transition and turbulent flow regimes in spite of the existence of the velocity potential. It was pointed out in [L-L] that transition to turbulence is related to the flow losing stability and this can be treated as a bifurcation process, i.e., period doubling which follows the Feigenbaum scenario [G]. Under a certain condition, the solution (4.1) becomes non-periodic ("chaotic") but bounded and nearby solutions separate rapidly in time. This latter property, called sensitive dependence upon the initial condition, can be thought as a loss of memory of the flow of its past history. It implies that long term predictions of the flow are almost impossible despite the deterministic nature of the Navier-Stokes equation.

We note further that the parameter $\Delta p/2\mu$ in the first integral of (4.1) is nothing else but $2\text{Re}/t$, where Re denotes the Reynolds number. Hence, as expected, Re is the parameter that quantifies the flow regime. We notice here that, because $2\text{Re}/t = 8\nu[\nabla \psi]^2/\psi$ , the pressure field can be fully excluded from the set of equations governing the velocity field. In other words, the wave function, $\psi$, if determined, fully defines the flow velocity.

Furthermore, as the first integral in (4.1) remains bounded for any $t$ and Re is finite, it follows that as $t \to \infty$, the first term in the right side of (4.1) becomes negligible in comparison with the second term. It means, in turn, that, provided a long enough time is given to a turbulent flow to exist, the flow will re-laminarize. This conclusion, drawn from observing (4.1), is consistent with the one drawn in [H et al], where it is claimed that the lifetimes of turbulence in pipe flow do not diverge at a finite critical Reynolds number. All the available data indicate that turbulence in pipe flows decays, and that eventually the flow will always re-laminarize. Similarly, lifetimes for the Couette flow and for a shear flow model do not appear to diverge, suggesting that a finite lifetime of turbulence may be a universal property of this class of flows.

A detailed analysis of (4.1) with respect to the phenomenon of laminar-turbulent transition will be presented in our next work.

## 5. Some general remarks on the solution

The present work establishes a fundamental link between the Schrödinger equation, which governs energy transport processes on quantum scales and the Navier-Stokes equation, which models phenomena of viscous fluid flow. It seems the very fact that the Navier-Stokes equation reduces to the Einstein-Kolmogorov equation (a non-homogeneous heat equation) is of fundamental importance, too. Indeed, if a fluid flow is seen from the topological perspective, it can be perceived as a space, whose curvature varies (in the case of turbulent flows, tremendously) with time and location. Yet, provided a sufficiently long time is given to the flow to exist, the solution, obtained here, predicts re-laminarization of the flow. In other words, in the course of time, smoothing of the flow space takes place; a direct consequence of the fact that the solution is written in terms of the Green's function of the heat equation. This, in turn, provides a striking similarity with the recent work by Perelman [P 1-3], who proved the Poincaré-Perelman theorem, in particular, using the fact that the evolution equation for the metric tensor implies the evolution equation for the curvature tensor, *Rm*, of the form

$$\frac{\partial Rm}{\partial t} = \nabla^2 Rm + Q$$

(5.1)

where $Q$ is a certain quadratic expression of the curvatures. In our opinion, a similarity between the approaches in [P 1-3] and this work is not a coincidence, but reflects the very fundamental properties of matter and the governing laws of nature.

# 6. Conclusions

The existence of the velocity potential is a direct consequence of the derivation of the continuity equation from the Schrödinger equation. This implies that the Cole-Hopf transformation is applicable to the Navier-Stokes equation for an incompressible flow and allows reducing the Navier-Stokes equation to the Einstein-Kolmogorov equation, in which the reaction term depends on the pressure. The solution to the resulting equation, and to the Navier-Stokes equation, can then be written in terms of the Green's function of the heat equation of the heat equation and is given in the form of an integral mapping. Such a form of the solution makes bifurcation period doubling possible, i.e. solutions for transition and turbulent flow regimes in spite of the existence of the velocity potential. One of the main properties of the solution, obtained in this work, is that, provided a long enough time is given to a turbulent flow to exist, the flow will re-laminarize.

## References:

[B]  F. J. Blatt, *Modern Physics* (1992) McGraw Hill, pp. 135 – 136.
[C]  J. D. Cole, Quart. Appl. Math. (1951) 9, 225.
[F]  U. Frisch, *Turbulence: the legacy of A. N. Kolmogorov* (1995) Cambridge University Press.
[G]  P. Glendinning, *Stability, Instability and Chaos: an Introduction to the Theory of Non-Linear Differential Equations* (1995) Cambridge University Press, Cambridge.
[H]  E. Hopf, Comm. Pure. Appl. Math. (1950) 3, 201.

[Hbrg]  W. Heisenberg, *The Physical Principles of the Quantum Theory* (1949) Dover.

[Hman] R. Haberman, *Elementary Applied Partial Differential Equations* (1987) Prentice Hall.

[H et al]  B. Hof et al, Finite lifetime of turbulence in shear flows (2006) Nature 443, pp. 59-62

[K-L]  V. V. Kulish, J. L. Lage, On the Relationship between Fluid Velocity and de Broglie's Wave Function and the Implications to the Navier-Stokes Equation (2002) Intl. J. Fluid Mech. Res. 29(1), pp. 40-52.

[L]  J. D. Logan, *An Introduction to Non-Linear Partial Differential Equations* (1994) John Wiley & Sons.

[Lgt]  A. J. Leggett, Superfluidity (1999) *Reviews of Modern Physics*, 71(2), pp. 318-323.

[L-L]  L. D. Landau, E. M. Lifshitz, *Fluid Mechanics*, 2nd Ed. (1987) Reed.

[P 1] G. Perelman, The entropy formula for the Ricci flow and its geometric applications. arXiv:math.DG/0211159 v1

[P 2] G. Perelman, Ricci flow with surgery on three-manifolds. arXiv:math.DG/0303109 v1

[P 3]  G. Perelman, Finite extinction time for the solutions to the Ricci flow on certain three-manifolds. arXiv:math.DG/0307245 v1

[T-S]  A. N. Tikhonov, A. A. Samarskii, Equations of Mathematical Physics (1963) Dover.